\documentclass[12pt,reqno]{amsart}
\usepackage{amsmath,bm}
\usepackage{amsfonts}
\usepackage{amssymb}
\usepackage{amsthm}
\usepackage{amscd}
\usepackage{a4wide}

\newcommand\R{\mathbb{R}}
\newcommand\Z{\mathbb{Z}}
\newcommand\N{\mathbb{N}}
\newcommand\Rp{\R^+}

\newcommand{\cal}[1]{\mathcal{#1}}

\newcommand{\cC}{{\cal C}}
\newcommand{\cD}{{\cal D}}

\newcommand{\cH}{{\cal H}}

\newcommand{\cS}{{\cal S}}

\newcommand{\vv}[1]{{\mathbf{#1}}}

\newtheorem{thm}{Theorem}

\newtheorem{cor}{Corollary}

\newtheorem{thkh}{Theorem\!\!}

\newtheorem{lemmabdv}{Lemma BDV\!\!}

\theoremstyle{remark}

\theoremstyle{definition}

\def\hs{{\mathcal H}^s}

\begin{document}

\title{ Inhomogeneous Diophantine approximation on \linebreak planar curves}

\author{V.V. Beresnevich, R.C. Vaughan}
\thanks{VVB: EPSRC Advanced Research Fellow: EP/C54076X/1}
\thanks{RCV: Research supported in part by NSA grant MDA904-03-1-0082 and  by a
Distinguished Visiting Fellowship at the University of York.}

\author{S.L. Velani}
\thanks{SLV: Research supported by EPSRC grants EP/E061613/1 and
EP/F027028/1 }

\begin{abstract}
The inhomogeneous metric theory for the set of simultaneously
$\psi$-approximable points lying on a planar curve is developed. Our
results naturally  incorporate  the  homogeneous
Khintchine-Jarn\'{\i}k type theorems recently established in
\cite{Beresnevich-Dickinson-Velani-07:MR2373145}
and~\cite{Vaughan-Velani-2007}.  The key lies in obtaining
essentially the best possible results regarding  the distribution of
`shifted' rational points near planar curves.
\\[2ex]
\end{abstract}

%
%

\date{\today}

\maketitle

%
%

\vspace*{3ex}

\section{{Introduction and Statement of Results} }

\subsection{Inhomogeneous approximation in the plane}

Throughout $\psi:\N\to\Rp$ is a monotonic function such that
$\psi(t)\to0$ as $t\to\infty$ and will be referred to as an
\emph{approximating function}. Given $\psi$ and a point
$\bm\theta:=(\theta_1,\theta_2) \in \R^2$, let $\cS(\psi,
\bm\theta)$ denote the set of points $\vv x :=(x_1,x_2) \in\R^2$ for
which there exists infinitely many positive integers $q $ such that
\begin{equation}\label{e:001}
\max_{1\leq i \leq2} \|qx_i-\theta_i\|<\psi(q)   \ .
\end{equation}
Here and throughout $\|\cdot\|$ denotes the distance to the
nearest integer. In the case that the inhomogeneous factor
$\bm\theta $ is the origin, the corresponding  set $\cS(\psi)$ is
the usual homogeneous set  of simultaneously $\psi$-approximable
points in the plane. In the case $\psi:t\to t^{-v}$ with $v>0$,
let us write $\cS(v, \bm\theta)$ for $\cS(\psi, \bm\theta)$. The
following statement provides a beautiful and simple criterion for
the `size' of $\cS(\psi, \bm\theta)$ expressed in terms of
$s$--dimensional Hausdorff measure $\hs$.

\begin{thkh}
Let $s \in (0,2]$, $\bm\theta \in \R^2$ and $\psi$ be  an
approximating function. Then
$$ \hs\left(\cS(\psi, \bm\theta)\right)=\left\{\begin{array}{ll} 0
& {\rm when} \;\;\; \sum \; t^{2-s} \,  \psi(t)^s \;\;
 <\infty\\ &
\\ \hs(\R^2) & {\rm when} \;\;\; \sum \;
 t^{2-s} \,  \psi(t)^s  \;\;  =\infty
\end{array}\right..$$
\end{thkh}

\noindent This result generalizes and unifies the classical
theorems of Khintchine (1924) and Jarn\'{\i}k (1931) and will be
refereed to as the Khintchine-Jarn\'{\i}k theorem. When $s=2$, the
measure ${\cal H}^2$ is equivalent to two-dimensional Lebesgue
measure in the plane and, loosely speaking, the theorem
corresponds to Khintchine's Theorem.
Actually,  the stronger statement that ${\cal H}^2(\R^2 \setminus
\cS(\psi, \bm\theta))=0$ if $\sum \psi(t)^2 = \infty$ is true and
the homogeneous case of this statement is due to  Khintchine. When
$s < 2$, the homogeneous case of the theorem corresponds to
Jarn\'{\i}k's Theorem and can be regarded as a Hausdorff measure
version of Khintchine's Theorem. For further details see
\cite[Section 12.1]{Beresnevich-Dickinson-Velani-06:MR2184760} and
references within.

\subsection{Inhomogeneous approximation restricted to curves} Let $\cC$ be a planar curve.
 In short, the goal is to obtain an analogue of
the above Khintchine-Jarn\'{\i}k theorem  for $\cC \cap \cS(\psi, \bm\theta) $. The fact
that the points $\vv x :=(x_1,x_2) \in\R^2$ of interest are
restricted to $\cC$ and therefore are of dependent
variables, introduces major difficulties in attempting to
describe
the measure theoretic structure of  $\cC \cap \cS(\psi, \bm\theta)$.

In 1998, Kleinbock $\&$ Margulis
\cite{Kleinbock-Margulis-98:MR1652916} established the fundamental
Baker-Sprindzuk conjecture concerning homogeneous Diophantine
approximation on  manifolds.  As a consequence, for non-degenerate
planar curves\footnote{A planar curve ${\cal C}$ is non-degenerate
if   the set of points on ${\cal C}$ at which the curvature
vanishes is a set of one--dimensional Lebesgue measure zero.
Moreover, it is not difficult to show that the set  of points on a
planar curve at which the curvature vanishes but the curve  is
non-degenerate is at most countable. In view of this, the
curvature completely describes the non-degeneracy of planar
curves.} the one--dimensional Lebesgue measure ${\cal H}^1$ of the
set $\cC \cap \cS(v)$ is zero  whenever $v > 1/2 $ -- see also \cite{schmidt}. Subsequently,
staying within the  homogeneous setup,  the significantly stronger
Khintchine-Jarn\'{\i}k type theorem for $\cC \cap \cS(\psi) $ has
been established -- see
\cite{Beresnevich-Dickinson-Velani-07:MR2373145} for the
convergence  part and \cite{Vaughan-Velani-2007} for the
divergence part.

Until the recent proof of the inhomogeneous Baker-Sprindzuk
conjecture
\cite{Beresnevich-Velani-Moscow,Beresnevich-Velani-08-Inhom}, the
theory of inhomogeneous Diophantine approximation on planar curves
(let alone manifolds) had remained essentially non-existent and
ad-hoc. As a consequence of the measure results in
\cite{Beresnevich-Velani-Moscow,Beresnevich-Velani-08-Inhom} or alternatively the even more recent dimension results in \cite{dmitry}, we
now know that for any non-degenerate planar curve $\cC$  and
$\bm\theta \in \R^2$,
$$ {\cal H}^1 (\cC \cap \cS(v, \bm\theta)) \, = \, 0 \qquad  {\rm when}  \quad  v > 1/2  \ . $$
Clearly, this statement is far from  the desirable
Khintchine-Jarn\'{\i}k type theorem for $\cC \cap \cS(\psi,
\bm\theta) $. As mentioned above, such a statement exists within
the homogeneous setup. This paper constitutes part of a programme
to develop a coherent inhomogeneous theory for curves, and indeed
manifolds in line with the homogeneous theory.

Without loss of generality,  we will assume that $\cC=\cC_f  : =
\{(x,f(x)):x\in I \} $ is given as the graph of a function
$f:I\to\R$, where $I$ is some interval of $\R$. As usual,
$C^{(n)}(I)$ will denote the set of $n$-times continuously
differentiable functions defined on some interval $I$ of $\R$.
In this paper we establish the inhomogeneous analogues of the main
theorems in \cite{Beresnevich-Dickinson-Velani-07:MR2373145} and
\cite{Vaughan-Velani-2007}; that is, we obtain the following
complete  Khintchine-Jarn\'{\i}k type theorem for planar curves.

\begin{thm}\label{thm1}
Let $s \in (1/2,1]$, $\bm\theta \in \R^2 $ and $\psi$ be an
approximating function. Let $f\in C^{(3)}(I)$ and assume that $
\hs(\{x\in I: f''(x)=0\}) = 0$.  Then
$$
\hs\left({\cC}_f\cap \cS(\psi, \bm\theta) \right)=\left\{
\begin{array}{ll} 0 & {\rm when} \;\;\; \sum
\; t^{1-s} \,  \psi(t)^{s+1} \;\;
 <\infty\\[1ex] &
\\ \cH^s(\cC_f) & {\rm when} \;\;\; \sum \;
 t^{1-s} \,  \psi(t)^{s+1}  \;\;  =\infty
\end{array}
\right..$$
\end{thm}

\medskip

\noindent   Note that a planar curve is one dimensional and so
$\hs\left({\cC}_f\cap \cS(\psi, \bm\theta) \right) \leq \hs({\cC}_f)
= 0 $ for any  $s> 1$ irrespective of the approximating function
$\psi$. Thus the hypothesis that  $ s \leq 1 $ is essential and
obvious.  In the case $s=1$,  the theorem is a statement concerning
the one-dimensional Lebesgue measure  of ${\cC}_f\cap \cS(\psi,
\bm\theta)$ and the convergence part actually only requires that
$f\in C^{(2)}(I)$. Also, as one would expect, the measure zero
assumption  on the set $\{x\in I: f''(x)=0\}$ coincides with the
definition of non-degeneracy. In the case $s<1$,  we have that
$\cH^s(\cC_f)=\infty$ and the theorem provides an elegant
zero-infinity law for the Hausdorff measure of
 ${\cC}_f\cap \cS(\psi, \bm\theta)$. In particular, this
law implies the following  corollary on the Hausdorff dimension of
${\cC}_f\cap \cS(\psi, \bm\theta)$ expressed in terms of the lower
order $ \lambda_\psi$ of $1/\psi$. Recall, $$ \lambda_\psi \, :=
\, \liminf_{t\to\infty}\frac{-\log \psi(t)}{\log t}
$$ and indicates the growth of the function $1/\psi$ `near'
infinity. Note that $\lambda_\psi$ is non-negative since
$\psi(t)\to0$ as $t\to\infty$.

\begin{cor}
Let $\bm\theta \in \R^2 $ and  $\psi$ be an approximating function
such that $\lambda_\psi\in [1/2,1) $. Let  $f \in C^{(3)}(I)$ such
that $f''(x)$ is not identically zero and assume that
$$
\dim \left\{  x \in I : f''(x) = 0  \right\}   \  \leq \
\frac{2-\lambda_\psi}{1+\lambda_\psi}  \ .
$$
Then,
$$ \dim{\cC}_f\cap \cS(\psi, \bm\theta)  = \frac{2-\lambda_\psi}{1+\lambda_\psi} .
$$
\end{cor}

\noindent This generalizes Theorem~4 of
\cite{Beresnevich-Dickinson-Velani-07:MR2373145} to the inhomogeneous
setting. Note that when $\lambda_{\psi}<1/2$, the condition
that $f''(x)$ is not identically zero follows from the assumption  that
$\dim \left\{  x
\in I : f''(x) = 0 \right\}   \,   < \,  1 $. We take this opportunity to mention that this necessary condition should also be present
in   \cite[Theorem~4]{Beresnevich-Dickinson-Velani-07:MR2373145},
where it is missing.

\vspace*{4ex}

\subsection{The  inhomogeneous counting results}

\noindent The proof of Theorem~\ref{thm1} rests on understanding  the
distribution of `shifted' rational points `near' planar  curves.
In view of the metrical nature of Theorem~\ref{thm1}, there is no harm is
assuming that the function $f :I\to\R$ is defined on a closed interval $I$
and that $f''$ is continuous and non-vanishing on $I$. By the compactness of $I$, there exist
positive and finite
constants $c_1,c_2$ such that
\begin{equation}\label{e:002}
 c_1 \ \le \  |f''(x)| \ \le \ c_2\qquad\forall \ x\in I.
\end{equation}

Let  $I$ and $f$ be as above.  Furthermore, given  $\bm\theta =(\theta_1,\theta_2)\in \R^2$,
$\delta > 0 $ and $Q\ge 1$, consider  the counting function
\begin{equation}\label{e:006}
N_{\bm\theta}(Q,\delta) \ := \ {\rm{card}}\left\{ (p_1,q) \in \Z \times {\mathbb N} \, : \,
\begin{array}{l}
 Q < q\le 2Q,\ (p_1+\theta_1)/q\in I\\[.5ex]
 \|qf(\,(p_1+\theta_1)/q\,)-\theta_2\|<\delta
 \end{array}
 \right\} \ .
\end{equation}
In short, the
function $N_{\bm\theta}(Q,\delta)$ counts the number of rational
points $(p_1/q,p_2/q)$ with bounded denominator $q$ such that the
shifted points $((p_1+\theta_1)/q,(p_2+\theta_2)/q)$  lie
within the $\delta/Q$-neighborhood of the
curve $\mathcal{C}_f$. The following result generalizes Theorem~1 of
\cite{Vaughan-Velani-2007} to the inhomogeneous setting.

\begin{thm}
\label{t:thm1}
Let $f\in C^{(2)}(I)$. Suppose that $Q\ge 1$ and $0<\delta\le\frac12$.  Then
\begin{equation*}
N_{\bm\theta}(Q,\delta) \ll \delta Q^2 + \delta^{-\frac12}Q  \ .
\end{equation*}
\end{thm}

With a mild additional condition on $f$ we are able to extend
the validity of the bound in Theorem \ref{t:thm1}. The following statement
is the inhomogeneous analogue of  Theorem~3 in
\cite{Vaughan-Velani-2007}.

\begin{thm}
\label{t:thm4}
Let $f''\in{\rm{Lip}}_{\phi}(I)$, where $0<\phi<1$.
Suppose that $Q\ge 1$  and
$0<\delta\le \frac12$.  Then, for any $\varepsilon > 0$
\begin{equation*}
N_{\bm\theta}(Q,\delta) \ll \delta Q^2 +
\delta^{-\frac12}Q^{\frac12+\varepsilon} +
\delta^{\frac{\phi-1}2}Q^{\frac{3-\phi}2}   \ .
\end{equation*}
\end{thm}

\medskip

\noindent {\em Remark.} When $\phi=1$ the proof gives the above
theorem with the term $\delta^{\frac{\phi-1}2}Q^{\frac{3-\phi}2}$
replaced by $Q\log(Q/\delta)$, and this is then always bounded by
one of the other two terms.

\medskip

 Armed with Theorems \ref{t:thm1} and  \ref{t:thm4}, the convergent
part of Theorem \ref{thm1} is established on following the
arguments set out in Sections 6 and 7 of
\cite{Vaughan-Velani-2007}. The modifications are essentially
obvious and the details are omitted. It is worth mentioning that
when $s=1$, we only need to appeal to Theorem \ref{t:thm1} and
thus we only require that  $f\in C^{(2)}(I)$ when proving the
convergent part of Theorem \ref{thm1}.

\medskip

  The key to establishing the divergence part
of Theorem~\ref{thm1} is the following covering result that also
yields  a  sharp lower bound for the counting function
$N_{\bm\theta}(Q,\delta)$. Throughout, $|X|$ will denote the
one-dimensional Lebesgue measure of a set $X$ in $\R$.

\begin{thm}\label{thm6}
Let $f\in C^{(3)}(I)$.
Then for any interval $J\subseteq I$ there are constants
$k_1,k_2,C_1,Q_0>0$ such that for any choice of $\delta$ and
$Q>Q_0$ subject to
\begin{equation}\label{e:007}
    \frac{k_1}{Q}\le\delta\le k_2
\end{equation}
one has
\begin{equation}\label{vb+}
\left|\bigcup_{(p_1,q)\in A_{\bm\theta}(Q,\delta,J)}
\left(B\left(\frac{p_1+\theta_1}{q},\frac{C_1}{Q^2
\delta}\right)\cap J\right)\right| \ \ge \  \frac{1}{2} \,
|J|\qquad\text{$\forall$ $\bm\theta\in\R^2$,}
\end{equation}
where
$$
A_{\bm\theta}(Q,\delta,J)\ :=\ \left\{(p_1,q)\in\Z\times \N \, : \,
\begin{array}{l}
Q<q\le 2Q, \ (p_1+\theta_1)/q\in J   \\[.5ex]
\|qf(\,(p_1+\theta_1)/q \,) - \theta_2 \|< \delta \end{array}
 \right\} \ \ .
$$
\end{thm}

\noindent Theorem \ref{thm6} is the inhomogeneous generalization of
Theorem~7 in \cite{Beresnevich-Dickinson-Velani-07:MR2373145}. Armed
with Theorem \ref{thm6},  the arguments set out in Section 7 of
\cite{Beresnevich-Dickinson-Velani-07:MR2373145} are easily adapted
to prove the divergence part of Theorem \ref{thm1}. The minor
modifications are essentially obvious and the details are omitted.


Note that  $N_{\bm\theta}(Q,\delta)$ is by definition the
cardinality of $A_{\bm\theta}(Q,\delta,I)$.  With this in mind, it
trivially  follows that
\begin{eqnarray*}
N_{\bm\theta}(Q,\delta) \cdot \frac{2C_1}{Q^2 \delta}  \  \  &
\geq & \ \sum_{ (p_1,q) \in A_{\bm\theta} (Q,\delta,I) } \left|
B\left(\frac{p_1+\theta_1}{q},\frac{C_1}{Q^2 \delta}\right)
\right|  \\[2ex]
&  \geq &  \left|\bigcup_{(p_1,q)\in A_{\bm\theta}(Q,\delta,I)}
\left(B\left(\frac{p_1+\theta_1}{q},\frac{C_1}{Q^2
\delta}\right)\cap I\right)\right|  \ \
\stackrel{(\ref{vb+})}{\ge} \  \ \frac{1}{2} \, |I| \ .
\end{eqnarray*}

\noindent In other words,  Theorem~\ref{thm6} implies the
following statement which is a generalisation of Theorem~6 in
\cite{Beresnevich-Dickinson-Velani-07:MR2373145}.

\begin{thm}\label{thm3}
Let $f\in C^{(3)}(I)$. There are constants $k_1,k_2,c,Q_0>0$ such
that for any choice of $\delta$ and $Q>Q_0$ satisfying
(\ref{e:007}) we have
$$
 N_{\bm\theta}(Q,\delta) \ \ge \  c \,  \delta\, Q^2 \qquad\text{$\forall$ $\bm\theta\in\R^2$}.
$$
\end{thm}

\vspace*{2ex}

\section{The proof of Theorems \ref{t:thm1} and \ref{t:thm4} }

\noindent Without loss of generality, assume that
$\bm\theta=(\theta_1,\theta_2)$ satisfies $0\le \theta_1, \theta_2
<1$. Let
\begin{equation*}\label{e:008}
J:=\left\lfloor \frac1{2\delta} \right\rfloor
\end{equation*}
and consider the Fej\'er kernel
\begin{equation*}
\mathcal K_J(x) := J^{-2}\left| \sum_{h=1}^Je(hx) \right|^2  =
\left( \frac{\sin\pi Jx}{J\sin\pi x} \right)^2.
\end{equation*}
When $\|x\|\le \delta$ we have $|\sin\pi Jx| = \sin\pi\|Jx\|\ge
2\| Jx\| = 2\|\, J\|x\|\, \|=2J\|x\|$, since $J\|x\| \le \delta
\left\lfloor \frac1{2\delta} \right\rfloor \le \frac12$.  Hence,
when $\|x\|\le\delta$, we have
\begin{equation*}
\mathcal K_J(x) \ge \frac{2\|x\|J}{J\pi\|x\|}=\frac2{\pi} \, .
\end{equation*}
Thus
\begin{equation*}
N_{\bm\theta}(Q,\delta)\le \frac{\pi}2 \sum_{Q<q\le 2Q}
\sum_{p_1+\theta_1\in qI} \mathcal K_J \big(
qf((p_1+\theta_1)/q)-\theta_2 \big).
\end{equation*}
Since
\begin{equation*}
\mathcal K_J(x) = \sum_{j=-J}^J \frac{J-|j|}{J^2} \, e(jx)
\end{equation*}
we have
\begin{equation*}
N_{\bm\theta}(Q,\delta)\le \pi\delta |I|Q^2 + N_1 +O(\delta Q) =
N_1 +O(\delta Q^2)
\end{equation*}
where $$N_1:= \frac{\pi}2 \sum_{0<|j|\le J} \frac{J-|j|}{J^2}
\sum_{Q<q\le 2Q} \  \sum_{p_1+\theta_1\in qI}
e\big(jqf((p_1+\theta_1)/q)-j\theta_2\big).$$ We observe that the
function $F(x):= jqf(x/q)$ has derivative $jf'(x/q)$. Given $j$
with $0<|j|\le J$ we define
$$H_-:=\lfloor\inf jf'(x)\rfloor-1,\ \
H_+:=\lceil\sup jf'(x)\rceil+1,$$
$$h_-:=\lceil\inf
jf'(x)\rceil+1,\ \ h_+:=\lfloor\sup jf'(x)\rfloor-1$$ where the
extrema are over $x$ in the interval $I$.  Then, by Lemma 4.2  of
\cite{Vaughan-97:MR1435742},
\begin{eqnarray*}
\sum_{p_1+\theta_1 \in qI} e\big(jqf((p_1+
\theta_1)/q)-j\theta_2\big) = \sum_{H_-\le h\le H_+}
  \int_{qI-\theta_1}
e\big(jqf((x \hspace*{-5ex} & & +\theta_1)/q)-j\theta_2-hx\big)dx \\
& + & O\big(\log(2+H)\big)
\end{eqnarray*}
where $H=\max(|H_-|,|H_+|)$.  Clearly $H\ll|j|\le J$ and so

$$N_1=N_2 + O\big(Q\log\textstyle{\frac1\delta}\big)$$
where
$$N_2
:= \frac{\pi}2 \sum_{0<|j|\le J} \frac{J-|j|}{J^2} \sum_{Q<q\le
2Q} \sum_{H_-\le h\le H_+} \int_{qI-\theta_1}
e\big(jqf((x+\theta_1)/q)-j\theta_2-h x\big)dx.$$

\noindent The integral here is
$$qe(h\theta_1-j\theta_2)\int_I e\big(q(jf(y)-hy)\big)dy.$$
As in Section 2  of \cite{Vaughan-Velani-2007}, we obtain
$$N_2=N_3+O\left(
\delta^{\frac12}Q^{\frac32} \right)$$ where
\begin{equation*}\label{e:009}
N_3:= \frac{\pi}2 \sum_{0<|j|\le J} \frac{J-|j|}{J^2} \sum_{Q<q\le
2Q} q \sum_{h_-< h < h_+} e(h\theta_1-j\theta_2) \int_I
e\big(q(jf(\theta_2)-h\theta_2)\big)d\theta_2
\end{equation*}
and the sum over $h$ is taken to be empty when $h_+\le h_-+1$.
Apart from the twisting factor $e(h\theta_1-j\theta_2)$ this
expression is identical to (2.3) of \cite{Vaughan-Velani-2007},
with identical properties of $f$.  The analysis of Sections 2 and
4 of \cite{Vaughan-Velani-2007}  can be applied without further
change to obtain the concomitant conclusions.


\section{The proof of Theorem~\ref{thm6}}

We will make use of the following result which appears as Lemma 6 in
\cite{Beresnevich-Dickinson-Velani-07:MR2373145}.

\begin{lemmabdv}\label{BKM0}
Let $\vv g :=(g_1,g_2):I\to\R^2$ be a $C^{(2)}$ map defined on a
compact interval $I$ such that $(g_1'g_2''-g_2'g_1'')(x)\neq0$  for
all $x\in I$. Given positive real numbers $\lambda,K,T$ and an
interval $J\subseteq I$, let $B(J,\lambda,K,T)$ denote the set  of
$x\in J$ for which there exists
$(q,p_1,p_2)\in\Z^3\smallsetminus\{0\}$ satisfying the following
system of inequalities:
$$
\left\{
\begin{array}{l}
|q \, g_1(x) \, + \, p_1 \, g_2(x)+ p_2| \ \le \ \lambda
 \\[2ex]
|q \, g_1'(x) \, + \, p_1 \, g_2'(x)| \ \le \ K
 \\[2ex]
|q| \ \le \  T \ \  .
\end{array}
\right.
$$
Then for any interval $J\subset I$ there is $C>0$ such that for any
choice of numbers $\lambda,K,T$ satisfying
\begin{equation}\label{e:010}
0<\lambda\le 1,\quad T\ge1, \quad K>0\quad \text{and} \quad
 \lambda KT\le1
\end{equation}
one has
\begin{equation}\label{e:011}
|{B(J,\lambda,K,T)}|\le C \max\left(\lambda^{1/3}, \left(\lambda K T
\right)^{1/9}\right)|J|   \ .
\end{equation}
\end{lemmabdv}

\vspace*{2ex}

To begin the proof of Theorem~\ref{thm6}, define  $\vv g(x) :=(g_1(x),g_2(x))$ by setting
$$
g_1(x):=x f'(x)-f(x)\qquad\text{and}\qquad g_2(x):=-f'(x).
$$
Then $\vv g\in C^{(2)}(I)$. Also, note that
\begin{equation}\label{e:012}
\vv g'(x)=(xf''(x),\,-f''(x)) \;  ,\qquad \vv
g''(x)=(f''(x)+xf'''(x),\,-f'''(x))
\end{equation}
and $$ (g_1'g_2''-g_2'g_1'')(x)=f''(x)^2 \ . $$ As $f''(x)\neq0$
everywhere, Lemma~BDV is applicable to this $\vv g$. In view
of  (\ref{e:002}) and the fact that $g_2'(x)=-f''(x)$,  it follows that
\begin{equation}\label{e:013}
 c_1\le |g_2'(x)|\le c_2  \qquad \forall x\in I \, .
\end{equation}
 For a fixed $x\in I$, consider the following
system of inequalities:
\begin{equation}\label{e:014}
\left\{
\begin{array}{l}
 |q^{}g_1(x)+p_1g_2(x)+p_2|\le c_0^3\delta
 \\[2ex]
 |q^{}g_1'(x)+p_1g_2'(x)|\le c_2(c_0^6Q\delta)^{-1}
 \\[2ex]
  |q|\le c_0^3Q \ .
\end{array}
\right.
\end{equation}
Here $c_0 < 1 $ is a real parameter to be determined later. Note that
with $q,p_1,p_2$ regarded as real variables,  the system defines a
convex body $\cD$ in $\R^3$ symmetric about the origin.

  Next, fix an interval $J \subseteq I $. By definition, the set $B(J,\lambda,K,T)$ with
\begin{equation}\label{e:015}
\lambda :=c_0^3\delta,\ \  \ K:=c_2(c_0^6Q\delta)^{-1},\ \ \
T:=c_0^4 Q
\end{equation}
consists of points $x\in J$ such that there exists a non-zero
integer solution $(q,p_1,p_2)$ to the system  (\ref{e:014}) with
$|q|\le c_0^4 Q$. By Lemma~BDV, for sufficiently large
$Q$ we have that
\begin{eqnarray*}
|B(J,\lambda,K,T)| & \le & C\, |J|\
\max\big\{(c_0^3\delta)^{1/3} \, , \ (c_2c_0)^{1/9}\big\} \\
& = & C\, (c_2c_0)^{1/9}|J| \ \le \ |J|/4 \  \
\end{eqnarray*}
provided that $c_0\le c_2^{-1}(4C)^{-9}$. Therefore,  with
$\lambda,K,T$ given by (\ref{e:015}) and $Q$ sufficiently large
\begin{equation}\label{e:016}
|{\textstyle\frac34}J\setminus B(J,\lambda,K,T)|\ge|J|/2  \ ,
\end{equation}
where $\frac34J$ is the interval $J$ scaled  by $\frac34$.

From this point onwards,   $x \in \frac34J\setminus B(J,\lambda,K,T)$ and is fixed. Then,
\begin{equation}\label{e:017} q>c_0^4 Q
\end{equation}
for any non-zero integer solution $(q,p_1,p_2)$ to the system  (\ref{e:014}). In other words,  the first consecutive minimum of the body (\ref{e:014})
is at least $c_0$. Let $\lambda_1\le\lambda_2\le\lambda_3$ be the
consecutive minima of the convex body $\cD$ given by (\ref{e:014}). Thus, $\lambda_1\ge
c_0$. By Minkowski's theorem on consecutive minima \cite{Cassels-97:MR1434478}, we have that
\begin{equation}\label{e:018}
\lambda_1\lambda_2\lambda_3V\le 2^3,
\end{equation}
where $V$ is the volume of $\cD$. It is readily
verified that
$$V \, =  \, 8|g_2'(x)|^{-1}c_2  \, \stackrel{(\ref{e:013})}{\ge}  \, 8c_2^{-1}c_2=8  \ . $$
Therefore,
$$
\lambda_3\stackrel{(\ref{e:018})}{\le} 8\lambda_1^{-2}V^{-1}\le
c_0^{-2}
$$

\noindent and it follows that there are three linearly independent
integer vectors
$$\vv a^{(i)}\, :=  \,  (q^{(i)},p_1^{(i)},p_2^{(i)})  \qquad   (1\le
i\le 3) $$ satisfying the system of inequalities
\begin{equation}\label{e:019}
\left\{
\begin{array}{l}
 |q^{(i)}g_1(x)+p_1^{(i)}g_2(x)+p_2^{(i)}|\le c_0\delta
 \\[2ex]
 |q^{(i)}g_1'(x)+p_1^{(i)}g_2'(x)|\le c_2(c_0^8Q\delta)^{-1}
 \\[2ex]
 0\le q^{(i)}\le c_0Q \ .
\end{array}
\right.
\end{equation}
For each $i$,  define
$$
G_i(x) \, : \, =q^{(i)}g_1(x)+p_1^{(i)}g_2(x)+p_2^{(i)}  \ .
$$
Now with $\bm\theta =(\theta_1,\theta_2)\in \R^2 $ fixed,    consider the following system of linear equations with respect to
the real variables $\eta_1,\eta_2,\eta_3$:
\begin{equation}\label{e:020}
    \left\{
     \begin{array}{rcl}
       \eta_1G_1(x)+\eta_2G_2(x)+\eta_3G_3(x) & = & \theta_1 f'(x)-\theta_2\\[2ex]
       \eta_1G'_1(x)+\eta_2G'_2(x)+\eta_3G'_3(x) & = & \theta_1 f''(x)\\[2ex]
       \eta_1q^{(1)}+\eta_2q^{(2)}+\eta_3q^{(3)} & = & 2Q    \ . \\[2ex]
     \end{array}
    \right.
\end{equation}
The determinant of this system is  equal to
$$
-f''(x)\left|\begin{array}{ccc}
               p_2^{(1)} & p_2^{(2)} & p_2^{(3)} \\[1ex]
               p_1^{(1)} & p_1^{(2)} & p_1^{(3)} \\[1ex]
               q^{(1)} & q^{(2)} & q^{(3)}
             \end{array}
\right|
$$
and so is non-zero. Therefore, there is a unique solution
$\eta_1,\eta_2,\eta_3$ to the system (\ref{e:020}). Let,
$t_i:=\lfloor\eta_i\rfloor$ when $q^{(i)}\ge0$ and
$t_i:=\lceil\eta_i\rceil$ when $q^{(i)}<0$. Therefore,
\begin{equation}\label{sv+}
|\eta_i-t_i|<1     \qquad   (1\le i\le 3)  \ .
\end{equation}
Also,  define $\vv
a=(q,p_1,p_2)\in\Z^3\setminus\{0\}$ by setting
\begin{equation}\label{e:021}
\vv a=(q,p_1,p_2):=\sum_{i=1}^3t_i\vv a^{(i)}.
\end{equation}

\noindent In view of the last equation of (\ref{e:020}) and the
definition of $t_i$, it follows  that $q\le 2Q$. Furthermore, using
the fact that $|q^{(i)}|\le c_0Q$ for each $i$ (this follows from
the last equation of (\ref{e:019})) we get that $q\ge 2Q-3c_0Q\ge Q$
provided that $c_0\le 1/3$. Thus, we have that
\begin{equation}\label{e:022}
    Q\le q\le 2Q.
\end{equation}
Further,
\begin{equation}\label{e:023}
\begin{array}[b]{rcl}
 |qg'_1(x)+p_1g'_2(x)-\theta_1 f''(x)| & \stackrel{(\ref{e:021})}{=} &
 \left|\sum_{i=1}^3t_iG'_i(x)-\theta_1 f''(x)\right|\\[2ex]
 & \stackrel{(\ref{e:020})}{=} &
 \left|\sum_{i=1}^3(t_i-\eta_i)G'_i(x)\right|\\[2ex]
 & \stackrel{(\ref{sv+})}{\le} &
 \sum_{i=1}^3|G'_i(x)|\\[2ex]
 & \stackrel{(\ref{e:019})}{\le} &
 3c_2(c_0^8Q\delta)^{-1}.
\end{array}
\end{equation}
In view of (\ref{e:012}) and (\ref{e:023}), we have that
$$
 |q x f''(x)-p_1f''(x)-\theta_1 f''(x)|<3c_2(c_0^8 Q\delta)^{-1}.
$$
The latter combined with (\ref{e:002}) gives that
\begin{equation}\label{e:024}
\left|x-\frac{p_1+\theta_1}{q}\right|\ \le \
\frac{3c_2}{q|f''(x)|c_0^8 Q\delta} \ \stackrel{(\ref{e:022})}{\le}
\ \frac{3c_2}{c_1c_0^8 Q^2\delta}=\frac{C_1}{Q^2\delta}\ ,
\end{equation}
where $C_1: =\frac{3c_2}{c_1c_0^8}$. For $Q$ sufficiently large, the
right hand side of (\ref{e:024}) can be made arbitrary small which
together with the fact that $x\in \frac34J$ ensures that
\begin{equation}\label{e:024+}
\textstyle\frac{p_1+\theta_1}{q}\in J  \ .
\end{equation}
\noindent Also,
\begin{equation}\label{e:025}
\begin{array}[b]{rcl}
 |q^{}g_1(x)+p_1g_2(x)+p_2-(\theta_1 f'(x)-\theta_2)| & \stackrel{(\ref{e:021})}{=} &
 \left|\sum_{i=1}^3t_iG_i(x)-(\theta_1 f'(x)-\theta_2)\right|\\[2ex]
 & \stackrel{(\ref{e:020})}{=} &
 \left|\sum_{i=1}^3(t_i-\eta_i)G_i(x)\right|\\[2ex]
 &\stackrel{(\ref{sv+})}{\le} &
 \sum_{i=1}^3|G_i(x)|\\[2ex]
 & \stackrel{(\ref{e:019})}{\le} &
 3c_0\delta.
\end{array}
\end{equation}
By Taylor's formula,
\begin{equation}\label{e:026}
\textstyle
f\big(\frac{p_1+\theta_1}{q}\big)=f(x)+f'(x)\big(\frac{p_1+\theta_1}{q}-x\big)+
\frac12f''(\tilde x)\big(\frac{p_1+\theta_1}{q}-x\big)^2
\end{equation}
for some $\tilde x$ between $x$ and $(p_1+\theta_1)/q$. Thus $\tilde x
\in J$. By (\ref{e:012}), the left hand side of (\ref{e:025}) equals
$|q(xf'(x)-f(x))-p_1f'(x)+p_2-(\theta_1 f'(x)-\theta_2)|$. Hence,
\begin{equation}\label{e:027}
\begin{array}[b]{rcl}
   3c_0\delta & \stackrel{(\ref{e:025})}{\ge}
& |q(xf'(x)-f(x))-p_1f'(x)+p_2-(\theta_1 f'(x)-\theta_2)| \\[2ex]
 &=&|(qx-p_1-\theta_1)f'(x)+p_2+\theta_2-qf(x)|\\[2ex]
 &\stackrel{(\ref{e:026})}{=}&
 \big|p_2+\theta_2-qf\big(\frac{p_1+\theta_1}{q}\big)+\frac{q}{2}f''(\tilde
x)(x-\frac{p_1+\theta_1}{q})^2\big| \, .
\end{array}
\end{equation}
Therefore, for $Q$ sufficiently large
\begin{equation}\label{e:028}
\begin{array}[b]{rcl}
 \big|qf\big(\frac{p_1+\theta_1}{q}\big)-p_2-\theta_2\big|&\le&
 \big| \, p_2+\theta_2-qf\big(\frac{p_1+\theta_1}{q}\big)
  + \frac{q}{2}f''(\tilde
 x)\big(x-\frac{p_1+\theta_1}{q}\big)^2\big|\\[2ex]
 &~ &  \hspace*{5ex} +  \ \  \big|\frac{q}{2}f''(\tilde
 x)\big(x-\frac{p_1+\theta_1}{q}\big)^2\big|\\[3ex]
 & \le & 3c_0\delta+\frac{Q}{2}c_2\left(\frac{C_1}{Q^2\delta}\right)^2
 \\[3ex]
 &\stackrel{(\ref{e:007})}{\le}&
 3c_0\delta+\frac{Q}{2}c_2\left(\frac{C_1}{k_1Q}\right)^2
 \, = \,  3c_0\delta+\frac{c_2C_1^2}{2k_1^2} \, Q^{-1}\\[3ex]
  &\stackrel{(\ref{e:007})}{\le}&
 3c_0\delta+\frac{c_2C_1^2}{2k_1^3} \, \delta \  \ <  \ \ \delta
\end{array}
\end{equation}
provided that $c_0<1/6$ and $\frac{c_2C_1^2}{2k_1^3}<1/2$. Thus, for
any $x \in \frac34J\setminus B(J,\lambda,K,T)$ there exists some
$(q,p_1,p_2)$ such that  (\ref{e:022}), (\ref{e:024+}) and
(\ref{e:028}) are satisfied.   In other words, $(q,p_1) \in
A_{\bm\theta}(Q,\delta,J)$ and moreover, in view of (\ref{e:016}) we
have that (\ref{vb+}) is satisfied for all $Q$ sufficiently large.
This completes the proof of Theorem~\ref{thm6}.

%
%
%

\vspace{1ex}

\noindent {\bf Acknowledgements.} SV would like to thank Resh Khodabocus for
his wonderful friendship and  expert support during the `broken
ankle' episode. Also, many thanks to the fab three -- Bridget,
Iona and Ayesha -- for putting up with an immobile and often bad
tempered man. Finally, SV would like to thank Victor Beresnevich
for his immense generosity (personally and professionally) and
Steve Donkin for being an extremely supportive chief during broken
times!

{\small \vspace{3mm}

\noindent Victor V. Beresnevich: Department of Mathematics,
University of York,

\vspace{-2mm}

\noindent\phantom{Victor V. Beresnevich: }Heslington, York, YO10
5DD, England.

\vspace{-1mm}

\noindent\phantom{Victor V. Beresnevich: }e-mail: vb8@york.ac.uk

\vspace{1mm}

\noindent Robert C. Vaughan: Department of Mathematics,
Pennsylvania State University

\vspace{-2mm}

\noindent\phantom{Robert C. Vaughan:  }University Park, PA 16802-6401,
U.S.A.

\vspace{-1mm}

\noindent\phantom{Robert C. Vaughan: }e-mail: rvaughan@math.psu.edu

\vspace{1mm}

\noindent Sanju L. Velani: Department of Mathematics, University of York,

\vspace{-2mm}

\noindent\phantom{Sanju L. Velani: }Heslington, York, YO10 5DD, England.

\vspace{-1mm}

\noindent\phantom{Sanju L. Velani: }e-mail:
slv3@york.ac.uk

}

\end{document}